\documentclass[12pt,a4paper]{amsart}

\usepackage[english]{babel}
\usepackage{amssymb}
\usepackage[leqno]{amsmath}
\usepackage{amsthm}
\usepackage{ae}
\usepackage[dvips]{graphicx}
\usepackage{amsmath}
\usepackage{comment}

\newcommand{\Z}{\mathbb{Z}}                  % The set of Integers.
\newcommand{\N}{\mathbb{N}}                  % Natural numbers.
\newcommand{\C}{\mathbb{C}}                  % Complex numbers.
\newcommand{\DD}{\mathbb{D}}                 % The unit disk
\newcommand{\R}{\mathbb{R}}          % Real line.
\newcommand{\Rn}{\mathbb{R}^n}               % Real n-space.
\renewcommand{\H}{{\bf H}}           % Halfspace
\newcommand{\p}{\partial}                    % Dho-symbol

\renewcommand{\a}{\alpha}
\renewcommand{\b}{\beta}
\newcommand{\g}{\gamma}
\renewcommand{\l}{\lambda}
\renewcommand{\d}{\delta}
\newcommand{\e}{\varepsilon}

\renewcommand{\o}{\omega}

\DeclareMathOperator{\dist}{dist}          % Distance.
\DeclareMathOperator*{\atan}{arc\,tan}     %   ``      ``    tangent.
\DeclareMathOperator{\Real}{Re}
\DeclareMathOperator{\Imaginary}{Im}

\def\Im{\Imaginary}
\def\Re{\Real}
\def\le{\leqslant}
\def\ge{\geqslant}

\def\hOmega{h_\Omega}%{(\Omega,h_\Omega)}
\def\kOmega{k_\Omega}%{(\Omega,k_\Omega)}
\def\hkOmega{\hork_\Omega}%{(\Omega,d_\Omega)}
\def\hork{\kappa}

\theoremstyle{plain}
\newtheorem{theorem}[equation]{Theorem}
\newtheorem{lemma}[equation]{Lemma}
\newtheorem{proposition}[equation]{Proposition}
\newtheorem{corollary}[equation]{Corollary}

\theoremstyle{definition}
\newtheorem{definition}[equation]{Definition}

\theoremstyle{remark}
\newtheorem{remark}[equation]{Remark}

\numberwithin{equation}{section}

\title[Gromov hyperbolicity of hyperbolic and
quasihyperbolic metrics]{Gromov hyperbolicity of Denjoy domains with hyperbolic and
quasihyperbolic metrics}

\author[P.\ H\"ast\"o]{Peter H\"ast\"o$^\ast$}
\address{Department of Mathematical Sciences,
P.O.\ Box 3000, FI-90014 University of Oulu, Finland}
\email{peter.hasto@helsinki.fi}
\thanks{$^\ast$ Supported in part by the Academy of Finland.}

\author[H.\ Lind\'en]{Henri Lind\'en}
\address{Department of Mathematics and Statistics, P.O.\ Box 64,
00140 University of Helsinki, Finland}
\email{hlinden@iki.fi}

\author[A.\ Portilla]{Ana Portilla$^\dagger$}
\address{St. Louis University (Madrid Campus),
Avenida del Valle 34, 28003 Madrid, Spain}
\email{apferrei@math.uc3m.es}
\thanks{$\dagger$ Supported in part by three grants from M.E.C.\
(MTM 2006-13000-C03-02, MTM 2006-11976 and MTM 2006-26627-E), Spain}

\author[J.\ M. Rodr{\'\i}guez]{Jos\'e M. Rodr{\'\i}guez$^{\;\dagger \ddagger}$}
\address{Departamento de Matem\'aticas, Universidad Carlos III de Madrid,
Avenida de la Universidad 30, 28911 Legan\'es, Madrid, Spain}
\email{jomaro@math.uc3m.es}
\thanks{$\ddagger$ Supported in part by a grant from
U.C.III$\,$M./C.A.M.\ (CCG06-UC3M/EST-0690), Spain.}

\author[E.\ Tour{\'\i}s]{Eva Tour{\'\i}s$^{\;\dagger \ddagger}$}
\address{Departamento de Matem\'aticas, Universidad Carlos III de Madrid,
Avenida de la Universidad 30, 28911 Legan\'es, Madrid, Spain}
\email{etouris@math.uc3m.es}

\date{\today}

\begin{document}


\begin{abstract} We obtain explicit and simple conditions
which in many cases allow one decide, whether or not a Denjoy domain
endowed with the Poincar\'e or quasihyperbolic metric is Gromov hyperbolic.
The criteria are based on the Euclidean size of the complement.
As a corollary, the main theorem
allows to deduce the non-hyperbolicity of any periodic Denjoy domain.
\end{abstract}

\maketitle

%%%%%% SECTION 1 %%%%%%%%%%%%%%%%%%%%%%

\section{Introduction}

In the 1980s Mikhail Gromov introduced a notion of abstract hyperbolic spaces,
which have thereafter been studied and developed by many authors. Initially,
the research was mainly centered on hyperbolic group theory, but lately
researchers have shown an increasing interest in more direct studies of
spaces endowed with metrics used in geometric function theory.

One of the primary questions is naturally whether a metric space $(X,d)$ is
hyperbolic in the sense of Gromov or not. The most classical examples,
mentioned in every textbook on this topic, are metric trees,
the classical Poincar\'e hyperbolic metric developed in the unit disk and,
more generally, simply connected complete
Riemannian manifolds with sectional curvature $K \le -k^2<0$.

However, it is not easy to determine whether a given space is Gromov
hyperbolic or not.
In recent years several investigators have been interested in showing that metrics
used in geometric function theory are Gromov hyperbolic.
For instance, the Klein-Hilbert metric
(see \cite{Be, KN}) is Gromov hyperbolic
(under particular conditions on the domain of definition);
that the Gehring-Osgood $j$-metric
(see \cite{Ha}) is Gromov hyperbolic;
and that the Vuorinen $j$-metric
(see \cite{Ha}) is not Gromov hyperbolic
except in the punctured space. Also, in \cite{Li} the hyperbolicity of
the conformal modulus metric $\mu$ and the related so-called Ferrand
metric $\lambda^*$, is studied.

Since the Poincar\'e metric is also the metric giving rise to what is commonly
known as the hyperbolic metric when speaking about open domains in the
complex plane or in Riemann surfaces, it could be expected
that there is a connection between the notions of hyperbolicity. For
simply connected subdomains $\Omega$ of the complex plane, it follows directly
from the Riemann mapping theorem that the metric space $(\Omega,h_\Omega)$
is in fact
Gromov hyperbolic. However, as soon as simple connectedness is omitted,
there is no immediate answer to whether the space $\hOmega$ is hyperbolic or
not. The question has lately been studied in
\cite{APRT} and \cite{PRT1}--\cite{RT3}.%,PRT2,PRT3,PorRTpp06,RT1,RT2,

The related quasihyperbolic metric has also recently been a topic of interest
regarding the question of Gromov hyperbolicity. In \cite{BHK}, Bonk, Heinonen and
Koskela found necessary and sufficient conditions for when a planar domain
$D$ endowed with the quasihyperbolic metric is Gromov hyperbolic.
This was extended by Balogh and Buckley, \cite{BB}:
they found two different necessary and sufficient conditions
which work in Euclidean spaces of all dimensions and also in
metric spaces under some conditions.

In this article we are interested in Denjoy domains. In this
case either the result of \cite{BHK} or \cite{BB} implies that
the domain is Gromov hyperbolic with respect to the
quasihyperbolic metric if and only if the domain is inner uniform
(see Section~\ref{positiveSect}). Although this is a vrey nice
characterization, it is somewhat difficult to check that a domain
is inner uniform, since we need to construct uniform paths
connecting every pair of points.

In this paper we show that it is necessary to look at paths joining
only a very small (countable) number of points when we want to
determine the Gromov hyperbolicity. This allows us to derive
a simple and very concrete conditions on when the domain is
Gromov hyperbolic. Much more importantly, our methods also suggest
corresponding results for the hyperbolic metric, which are also
proven. To the best of our knowledge, this is the first time
that Gromov hyperbolicity of any class of infinitely connected
domains has been obtained from conditions on the Euclidean
size of the complement of the domain.

The main results in this article are the following:

\begin{theorem}
\label{mainThm:qh}
\label{mainThm:h}
Let $\Omega$ be a Denjoy domain with $\Omega\cap\R =
(-\infty,0) \cup \bigcup_{n=1}^\infty (a_n,b_n)$, $b_n \le a_{n+1}$ for every
$n$, and $\lim_{n\to\infty} a_n = \infty$.
\begin{enumerate}
\item
The metrics $\kOmega$ and $\hOmega$ are Gromov hyperbolic if
\[
\liminf_{n\to\infty} \frac{b_n-a_n}{a_n} > 0. \]
\item
The metrics $\kOmega$ and $\hOmega$ are not Gromov hyperbolic if
\[\lim_{n\to\infty} \frac{b_n-a_n}{a_n} = 0. \]
\end{enumerate}
\end{theorem}

% \begin{theorem}
% \label{mainThm:h}
% Let us consider a Denjoy domain $\Omega$ with $\Omega\cap\R =
% (-\infty,0) \cup \bigcup_{n=1}^\infty (a_n,b_n)$, $b_n \le a_{n+1}$ for every
% $n$, and $\lim_{n\to\infty} a_n = \infty$.
% \begin{enumerate}
% \item
% If $\liminf_{n\to\infty} (b_n-a_n)/a_n > 0$, then
% $\hOmega$ is Gromov hyperbolic.
% \item
% If $\lim_{n\to\infty} (b_n-a_n)/a_n = 0$, then
% $\hOmega$ is not Gromov hyperbolic.
% \end{enumerate}
% \end{theorem}

It is interesting to note that
in the case of Denjoy domains many of the results seem to
hold for both the hyperbolic and the quasihyperbolic metrics.
In fact, we know of no planar domain which is Gromov hyperbolic with
respect to one of these metrics, but not the other.

In the previous theorems, the boundary components had a single accumulation
point, at $\infty$, and the accumulation happened only from one side.
It turns out that if this kind of domain is not Gromov hyperbolic,
then we cannot mend the situation by adding some boundary to
the other side of the accumulation point, as the following
theorem shows.

\begin{theorem}
\label{th:negativehalfaxis}
Let $\Omega$ be a Denjoy domain with $(-\infty,0)\subset \Omega$
and let $F\subseteq (-\infty,0]$ be closed.
If $\kOmega$ is not Gromov hyperbolic, then neither is
$k_{\Omega\setminus F}$; if $\hOmega$ is not Gromov hyperbolic, then neither is
$h_{\Omega\setminus F}$.
\end{theorem}

We also prove the non-hyperbolicity of any periodic Denjoy domain:

\begin{corollary}
\label{periodicCor}
Let $E_0 \subset [0,t)$ be closed, $t>0$, set $E_n:= E_0+tn$ for
$n\in \N$ or $n\in \Z$, and $\Omega:=\C \setminus \cup_{n} E_n$.
Then
$\hOmega$ and $\kOmega$ are not Gromov hyperbolic.
\end{corollary}

%%%%%%%%%%%%%%%%%%%%%%%%%%%%%%%%%%%%%%%%%%%%%%%%%%%%%%%%%%%%

\section{Definitions and notation}

By $\H^2$ we denote the upper half plane,
$\{z\in\C:\; \Imaginary z>0\}$, by $\DD$ the unit disk $\{z\in\C:\; |z|<1\}$.
For $D\subset \C$ we denote by $\partial D$ and $\overline{D}$ its
boundary and closure, respectively. For $z\in D\subsetneq \C$
we denote by $\delta_D(z)$ the distance to the boundary,
$\min_{a\in\partial D} |z-a|$.
Finally, we denote by $c$, $C$, $c_j$ and $C_j$ generic constants
which can change their value from line to line and even in the same line.

Recall that a domain $\Omega \subset \C$ is said to be of
\textit{hyperbolic type} if it has at least two finite boundary
points. The universal cover of such domain is the unit disk
$\DD$. In $\Omega$ we can define the Poincar\'e
metric, i.e.\ the metric obtained by pulling back the metric
$ds =2 |dz|/(1-|z|^2)$ of the unit disk.
Equivalently, we can pull back the metric $ds= |dz|/\Imaginary z$
of the the upper half plane $\H^2$. Therefore, any simply connected subset
of $\Omega$ is isometric to a subset of $\DD$. With this metric, $\Omega$ is
a geodesically complete Riemannian manifold with constant curvature
$-1$, in particular, $\Omega$ is a geodesic metric space. The Poincar\'e
metric is natural and useful in complex analysis; for instance, any
holomorphic function between two domains is Lipschitz with constant
$1$, when we consider the respective Poincar\'e metrics.

The quasihyperbolic metric is the distance induced by the density
$1/\delta_\Omega(z)$. By $\lambda_\Omega$ we denote the density of
the Poincar\'e metric in $\Omega$,
and by  $k_\Omega$ and $h_\Omega$ the quasihyperbolic and Poincar\'e distance
in $\Omega$, respectively. Length (of a curve) will be denoted by the
symbol $\ell_{d,\Omega}$, where $d$ is the metric with respect to which length
is measured. If it is clear which metric or domain is used, either one or both
subscripts in $\ell_{d,\Omega}$ might be left out. The subscript
$\text{Eucl}$ is used to denote the length with respect to the
Euclidean metric. Also, as most of the proofs
apply to both the quasihyperbolic and the Poincar\'e metrics, we will use
the symbol $\hork$ also as a ``dummy metric'' symbol, where it can be replaced
by either $k$ or $h$.

We denote by $\lambda_\Omega$ the density of the hyperbolic
metric in $\Omega$. It is well known that for every domain $\Omega$
$$
\l_\Omega (z) \le \frac2{\delta_\Omega (z)} \quad \;\forall\, z \in \Omega ,
\quad\qquad \ell_{h,\Omega}(\g) \le 2 \ell_{k,\Omega} (\g)
\quad \;\forall \, \g \subset \Omega ,
$$
and that for all domains $\Omega_1\subset \Omega_2$ we have
$\l_{\Omega_1} (z) \ge \l_{\Omega_2} (z)$ for every $z \in \Omega_1$.

If $\Omega_0$ is an open subset of $\Omega$, in $\Omega_0$ we always consider
its usual quasihyperbolic or Poincar\'e metric (independent of $\Omega$). If
$D$ is a closed subset of $\Omega$, we always consider in $D$
the inner metric obtained by the restriction of the
quasihyperbolic or Poincar\'e metric in $\Omega$, that is
\[
\begin{split}
d_{\Omega | D}(z,w):= &
\inf\big\{ \ell_{\hork,\Omega}(\g) : \, \g\subset D \text{ is a continuous }
\\
& \qquad \text{ curve joining $z$ and $w$} \big\}\ge d_{\Omega}(z,w)\,.
\end{split}
\]
It is clear that $\ell_{\Omega | D}(\g)=\ell_{\Omega}(\g)$
for every curve $\g\subset D$.
We always require that $\p D$ is a union of pairwise disjoint
Lipschitz curves; this fact guarantees that $(D, d_{\Omega | D})$ is a geodesic
metric space.

A geodesic metric space $(X,d)$ is said to be
{\it Gromov $\delta$-hyperbolic}, if
\[d(w,[x,z] \cup [z,y]) \le \delta\]
for all $x,y,z \in X$; corresponding
geodesic segments $[x,y],[y,z]$ and $[x,z]$; and $w \in [x,y]$.
If this inequality holds, we also say that the geodesic triangle
is {\it $\delta$-thin}, so Gromov hyperbolicity can be reformulated by requiring that
all geodesic triangles are thin.
% In non-geodesic spaces Gromov hyperbolicity
% can be defined by using the inequality
% $$(x|z)_w \ge \min\{(x|y)_w,(y|z)_w \}-\delta,$$
% where $w \in X$ is an arbitrary base point, and
% $$(x|y)_w=\frac{1}{2}\left(d(x,w)+d(y,w)-d(x,y) \right)$$
% denotes the {\it Gromov product} of $x$ and $y$. In this article we will only
% need the traditional geodesic wiewpoint.

A \textit{Denjoy domain $\Omega\subset\C$} is a domain whose boundary
is contained in the real axis. Hence, it satisfies
$\Omega\cap\R = \cup_{n\in \Lambda} (a_n,b_n)$,
where $\Lambda$ is a countable index set,
$\{(a_n,b_n)\}_{n\in \Lambda}$ are pairwise disjoint,
and it is possible to have
$a_{n_1}= -\infty$ for some $n_1 \in \Lambda$
and/or $b_{n_2}= \infty$ for some $n_2 \in \Lambda$.

In order to study Gromov hyperbolicity, we consider the case where
$\Lambda$ is countably infinite, since if $\Lambda$ is finite then $\hOmega$
and $\kOmega$ are easily seen to be Gromov
hyperbolic by Proposition~\ref{finitetype}, below.

%%%%%%%%%%%%%%%%%%%%%%%%%%%%%%%%%%%%%%%%%%%%%%%%%%%%%%%%%%%%
%%%%%%%%%%%%%%%%%%%%%%%%%%%%%%%%%%%%%%%%%%%%%%%%%%%%%%%%%%%%
%%%%%%%%%%%%%%%%%%%%%%%%%%%%%%%%%%%%%%%%%%%%%%%%%%%%%%%%%%%%

\section{Some classes of Denjoy domains which are Gromov hyperbolic}
\label{positiveSect}

The quasihyperbolic metric is traditionally defined in subdomains of
Euclidean $n$-space $\Rn$, i.e.\ open and connected subsets $\Omega \subsetneq
\Rn$. However, a more abstract setting is also possible, as was shown in
the article \cite{BHK} by Bonk, Heinonen and Koskela. There it is shown
that if $(X,d)$ is taken to be any metric space which is locally compact,
rectifiably connected and noncomplete, the quasihyperbolic metric $k_{X}$
can be defined as usual, using the weight $1/\dist(x,\p X)$.

%\begin{definition}
Given a real number $A \ge 1$, a curve $\gamma \colon [0,1] \to \Omega$ is called
 {\it $A$-uniform for the metric $d$} if
\begin{eqnarray*}
&{}& \ell_d(\gamma) \le A\; d(\g(0),\g(1)) \quad\text{and}\\
&{}& \min\{ \ell_d(\gamma|[0,t]),\ell_d(\gamma|[t,1])\} \le A\; \dist_d(\g(t),\p \Omega),\qquad {\rm for\; all}\; t \in [0,1].
\label{unif2}
\end{eqnarray*}
Moreover, a locally compact, rectifiably connected noncomplete metric space
is said to be {\it $A$-uniform} if every pair of points can be joined by
an $A$-uniform curve.
The abbreviations ``$A$-uniform'' and ``$A$-inner uniform''
(without mention of the metric) mean $A$-uniform for the
Euclidean metric and Euclidean inner metric, respectively.
%\end{definition}

Uniform domains are intimately connected to domains which are Gromov
hyperbolic with respect to the quasihyperbolic metric
(see \cite[Theorems~1.12, 11.3]{BHK}). Specifically, for a Denjoy domain $\Omega$
these results imply that $\kOmega$ is Gromov hyperbolic if and only
if $\Omega$ is inner uniform.

Here we will use the generalized setting in \cite{BHK} to show
that for Denjoy domains it actually suffices to consider the upper (or lower)
intersection with the actual domain, as can be done for the Poincar\'e metric:

\begin{lemma}
\label{l:halfplane}
Let $\emptyset \ne E \subset \mathbb{R}$ be a closed set, and denote
$D_0=\C \setminus E$ and $D=D_0 \cap \{z \in \C\; |\;
\Imaginary z \ge 0\} = D_0 \cap \overline{{\bf H}^2}$. Then the metric space
$D$, with the restriction of the Poincar\'e or the quasihyperbolic metric in $D_0$,
is $\d$-Gromov hyperbolic,
with some universal constant $\d$.
\end{lemma}

\begin{proof}
We deal first with the quasihyperbolic metric.
As the upper half-plane is uniform in the classical case, the same curve of
uniformity (which is an arc of a circle orthogonal to $\R$) can be shown to
be an $A$-uniform curve in the sense of \cite{BHK} for the set $D$. Hence $D$
is $A$-uniform. By \cite[Theorem 3.6]{BHK} it then follows that the space
$(D,k_D)$ is Gromov hyperbolic.

We also have that $D$ is hyperbolic with the restriction of the Poincar\'e
metric $h_{D_0}$, since it is isometric to a geodesically convex subset of the
unit disk (in fact, there is just one geodesic in $D$ joining two points in
$D$). Therefore, $D$ has $\log\big(1+\sqrt{2}\,\big)$-thin triangles, as does
the unit disk (see, e.g.\ \cite[p.~130]{An}).
\end{proof}

\begin{definition}
Let $\Omega$ be a Denjoy domain. Then we have
$\Omega \cap \R = \cup_{n\ge 0} (a_n,b_n)$ for some suitable intervals.
We say that a curve in $\Omega$ is a \emph{fundamental
geodesic} if it is a geodesic joining
$(a_0,b_0)$ and $(a_n,b_n)$, $n>0$,
which is contained in the closed halfplane
$\overline{\H^2}=\{z\in\C:\,\Imaginary z \ge 0\}$.
We denote by $\g_n$ a fundamental geodesic corresponding to $n$.
\end{definition}

The next result was proven for the hyperbolic metric in \cite[Theorem~5.1]{APRT}.
In view of Lemma~\ref{l:halfplane} one can check that the same proof carries
over to the quasihyperbolic metric.

By a \textit{bigon} we mean a closed polygon with two edges.

\begin{theorem}
\label{th:2.1}
Let $\Omega$ be a Denjoy domain and denote by $\hork_\Omega$
the Poincar\'e or quasihyperbolic metric.
Then the following conditions are equivalent:

\begin{enumerate}
\item[$(1)$] $\hork_\Omega$ is $\d$-hyperbolic.

\item[$(2)$] There exists a constant $c_1$ such that for every choice of
fundamental geodesics $\{\g_n\}_{n=1}^\infty$ we have $\hork_\Omega(z,\R) \le c_1$ for
every $z\in \cup_{n\ge 1} \g_n$.

\item[$(3)$] There exists a constant $c_2$ such that for a fixed choice of
fundamental geodesics $\{\g_n\}_{n=1}^\infty$ we have $\hork_\Omega(z,\R) \le c_2$ for
every $z\in \cup_{n\ge 1} \g_n$.

\item[$(4)$] There exists a constant $c_3$ such that every geodesic bigon
in $\Omega$ with vertices in $\R$ is $c_3$-thin.
\end{enumerate}

Furthermore, the constants in each condition
only depend on the constants appearing in any other of the conditions.
\end{theorem}

Note that the case $\Omega\cap\R = \cup_{n=0}^N (a_n,b_n)$
is also covered by the theorem.

\begin{corollary}
\label{c:2.2}
Let $\Omega$ be a Denjoy domain and denote by $\hork_\Omega$
the Poincar\'e or quasihyperbolic metric.
If there exist a constant $C$ and a sequence of fundamental geodesics
$\{\g_n\}_{n\ge 1}$ with $\ell_{\hork,\Omega}(\g_n) \le C$,
then $\hork_\Omega$ is $\d$-Gromov
hyperbolic, and $\d$ just depends on $C$.
\end{corollary}

% \begin{corollary}
% \label{c:finite}
% Let $\Omega$ be a Denjoy domain, denote by $d_\Omega$
% the Poincar\'e or quasihyperbolic metric and
% fix $N>0$. Then $\dOmega$ is Gromov hyperbolic
% if and only if
% there exists a constant $c$ such that for a fixed choice of
% fundamental geodesics $\{\g_n\}_{n\ge 1}$ we have $d_\Omega(z,\R) \le c$ for
% every $z\in \cup_{n> N} \g_n$.
% \end{corollary}

% \begin{proof}
% It is sufficient to apply the equivalence of $(1)$ and $(3)$ in Theorem
% \ref{th:2.1} with
% \[
% c_2:= \max \big\{c, \ell_{d,\Omega}(\g_1),\dots, \ell_{d,\Omega}(\g_N) \big\}. \qedhere
% \]
% \end{proof}

If $\Omega$ has only finitely many boundary components,
then it is always Gromov hyperbolic, in a quantitative way:

\begin{proposition}
\label{finitetype}
Let $\Omega$ be a Denjoy domain with
$\Omega\cap\R = \cup_{n=1}^N (a_n,b_n)$,
and denote by $\hork_\Omega$
the Poincar\'e or quasihyperbolic metric.
Then $\hork_\Omega$ is $\d$-Gromov hyperbolic,
where $\d$ is a constant which only depends on $N$
and $c_0 = \sup_n \hork_\Omega\big( (a_n,b_n), (a_{n+1},b_{n+1})\big)$.
\end{proposition}

Note that we do not require $b_n\le a_{n+1}$.

\begin{proof}
Let us consider the shortest geodesics $g_n^*$ joining
$(a_n,b_n)$ and $(a_{n+1},b_{n+1})$ in
$\Omega^+:= \Omega \cap \overline{\H^2}$.
Then $\ell_\Omega(g_n^*) \le \ell_\Omega(g_n) \le c_0$
for $0\le n \le N-1$.

By Theorem~\ref{th:2.1}, we just need to prove that there exists a
constant $c$, which only depends on $c_0$ and $N$, such that
$\hork_\Omega(z,\R) \le c$ for every $z \in \cup_{n= 1}^N \g_n$.

For each $0\le n \le N-1$, let us
consider the geodesic polygon $P$ in $\Omega^+$, with the following sides:
$\g_n,g_{0}^*,\dots,g_{n-1}^*$, and the geodesics
joining their endpoints which are contained in $(a_0,b_0),\dots,(a_{n},b_{n})$.
Since $(\Omega^+, \hork_\Omega)$ is $\d_0$-Gromov hyperbolic, where $\d_0$ is a constant
which only depends on $c_0$, by Lemma \ref{l:halfplane}, and $P$ is a
geodesic polygon in $\Omega^+$ with at most $2N+2$ sides, $P$ is $2N\d_0$-thin.
Therefore, given any $z\in \g_n$, there exists a point $w\in \cup_{k=0}^{N-1} g_k^*
\cup \R$ with $\hork_{\Omega}(z,w) \le 2N\d_0$. Since $\ell_\Omega(g_k^*)
\le c_0$ for $0\le k \le N-1$, there exists $x\in \R$ with
$\hork_{\Omega}(x,w) \le c_0/2$. Hence, $\hork_{\Omega}(z,\R) \le \hork_{\Omega}(z,x) \le 2N\d_0 + c_0/2$,
and we conclude that $\hkOmega$ is $\d$-Gromov hyperbolic.
\end{proof}

\begin{theorem}
\label{th:boundedlength}
Let $\Omega$ be a Denjoy domain with
$\Omega\cap\R = \cup_{n=0}^\infty (a_n,b_n)$,
$(a_0,b_0)=(-\infty,0)$ and $b_n \le a_{n+1}$ for every $n$.
Suppose that $b_n \ge K a_n$ for a fixed $K>1$ and every $n$.
Then $\hOmega$ and $\kOmega$ are $\d$-Gromov hyperbolic,
with $\d$ depending only on $K$.
\end{theorem}

\begin{proof}
Fix $n$ and consider the domain
$$
\Omega_n=\frac{1}{a_n}\, \Omega
=\left\{ \frac{x}{a_n}\; |\; x \in \Omega\right\}.
$$
If we define $D:= \C\setminus [0,1] \cup [K, \infty)$,
then $D \subset \Omega_n$, and
$\ell_{k,{\Omega_n}}(\gamma) \le \ell_{k,D}(\gamma)$
for every curve $\gamma \subset \Omega_n$.
The circle $\sigma:=S^1(0, (1+K)/2)$ goes around
the boundary component $[0,1]$ in $D$ and has finite quasihyperbolic
length:
$$
\ell_{k,D}(\sigma) \le \int_\sigma \frac{|dz|}{(K-1)/2}
= 2 \pi \, \frac{K+1}{K-1} \,.
$$
Consider the shortest fundamental geodesics
joining $(a_0,b_0)$ with $(a_n,b_n)$,
with the Poincar\'e and the
quasihyperbolic metrics, $\g_n^h$ and $\g_n^k$, respectively. Then,
$$
\begin{aligned}
\ell_{k,\Omega}(\g_n^k) &
= \ell_{k,{\Omega_n}}\Big(\frac{1}{a_n}\g_n^k\Big)
\le  \ell_{k,{\Omega_n}}(\sigma)
\le \ell_{k,D}(\sigma) \le  2 \pi \, \frac{K+1}{K-1} \,,
\\
\ell_{h,\Omega}(\g_n^h) & \le \ell_{h,\Omega}(\g_n^k) \le 2\,\ell_{k,\Omega}(\g_n^k)
\le 4 \pi \,
\frac{K+1}{K-1} \,.
\end{aligned}
$$
Therefore $\hOmega$ and $\kOmega$ are $\d$-Gromov
hyperbolic (and $\d$ depends only on $K$), by Corollary~\ref{c:2.2}.
\end{proof}

\begin{proof}[Proof of Theorems~\ref{mainThm:qh}(1)]
If $\liminf_{n\to\infty} (b_n-a_n)/a_n > 0$, then we can choose
$K>1$ so that $(b_n-a_n)/a_n > K-1$ for every $n$,
whence $b_n > K a_n$. Thus the previous theorem implies
the claims.
\end{proof}

%%%%%%%%%%%%%%%%%%%%%%%%%%%%%%%%%%%%%%%%%%%%%%%%%%%%%%%%%%%%%%%%
%%%%%%%%%%%%%%%%%%%%%%%%%%%%%%%%%%%%%%%%%%%%%%%%%%%%%%%%%%%%%%%%
%%%%%%%%%%%%%%%%%%%%%%%%%%%%%%%%%%%%%%%%%%%%%%%%%%%%%%%%%%%%%%%%

\section{Some classes of Denjoy domains which are not Gromov hyperbolic}

The following function was introduced by Beardon and Pommerenke \cite{BP}.

\begin{definition}
For $\Omega\subsetneq\C$, define $\beta_\Omega(z)$ as the function
$$
\beta_\Omega (z) := \inf\Big\{ \Big|\log\Big|\frac {z-a}{b-a}\Big|\Big|:
\ a,b\in\p\Omega, \ |z-a|=\delta_\Omega (z) \Big\} \, .
$$
\end{definition}

The function $\beta_\Omega$ has a geometric interpretation.
We say that an annulus $\{z\in\C: \, r< |z-a|<R \}$ \emph{separates}
$E\subset\C$ if
$\{z\in\C: \, r< |z-a|<R \} \cap E = \emptyset$, $\{z\in\C:
\, |z-a|\le r \} \cap E \neq \emptyset$ and $\{z\in\C: \,
|z-a|\ge R \} \cap E \neq \emptyset$.
We say that $E$ is \emph{uniformly perfect} if there exists a
constant $c_1$ such that $R/r \le c_1$ for every annulus
$\{z\in\C: \, r< |z-a|<R \}$ separating $E$ (see \cite{BP,P1,P2}).
Now we see that $\beta_\Omega$ is
bounded precisely when $\Omega$ is uniformly perfect.

Thus it follows from the next theorem, that
$\l_\Omega$ and $1/\delta_\Omega$ are comparable if and only if $\Omega$
is uniformly perfect.

\begin{theorem}[Theorem~1, \cite{BP}]\label{BPthm}
For every domain $\Omega\subset\C$ of hyperbolic type and
for every $z\in\Omega$, we have that
$$
2^{-3/2} \le \l_\Omega (z) \, \delta_\Omega (z) \, (k_0+\beta_\Omega(z)) \le \pi/4 \,,
$$
where $k_0=4+\log(3+2\sqrt{2})$.
\end{theorem}

\begin{lemma}\label{minLenLem}
Let $\gamma$ be a curve in a domain $D\subset \Rn$ from $a\in D$ with Euclidean
length $s$. Then:
\begin{enumerate}
\item
$\ell_{k,D}(\gamma) \ge \log \big( 1 + \frac{s}{d_D(a)} \big)\,$.
\item
If $D$ is a Denjoy domain and $a \in (a_n,b_n)$, with $b_n-a_n\le
r$, then $\ell_{h,D}(\gamma) \ge 2^{-3/2}\log \big( 1 + k_0^{-1} \log \big( 1 +
\frac{s}{r} \big) \big)$, with $k_0$ as in Theorem~\ref{BPthm}.
\end{enumerate}
\end{lemma}

\begin{proof}
Let $z\in \partial D$ be a point with $\delta_D(a)=|a-z|$. Without loss
of generality we assume that $z=0$. By monotonicity $\ell_{k,D}(\gamma) \ge
\ell_{k,{\Rn\setminus\{0\}}}(\gamma)$. Further, it is clear that
$\ell_{k,{\Rn\setminus\{0\}}}(\gamma)\ge \ell_{k,{\Rn\setminus\{0\}}}([|a|,
|a|+s])$, whence the first estimate by integrating the density $1/|x|$.

We then prove the second estimate.
Without loss of generality we assume that $b_n=0$.
By monotonicity $\ell_{h,D}(\gamma) \ge \ell_{h,{\C\setminus\{a_n,0\}}}
(\gamma)$.
By \cite[Theorem 4.1(ii)]{M} we have that
$\lambda_{\C\setminus\{a_n,0\}}(z) \ge \lambda_{\C\setminus\{a_n,0\}}(|z|)$
and by \cite[Theorem 4.1(i)]{M}
that $\lambda_{\C\setminus\{a_n,0\}}(r)$
is a decreasing function in $r\in(0,\infty)$; hence,
$\ell_{h,{\C\setminus\{a_n,0\}}}(\gamma) \ge \ell_{h,{\C\setminus\{a_n,0\}}}([|a_n|,|a_n|+s])
= \ell_{h,{\C\setminus\{-1,0\}}}([1,1+s/|a_n|])$.
By Theorem~\ref{BPthm}
\[
\begin{split}
\ell_{h,D}(\gamma)
& \ge \ell_{h,{\C\setminus\{-1,0\}}}([1,1+s/|a_n|])
\ge \int_{1}^{1+s/|a_n|} \frac{2^{-3/2} \, dx}{x\,\big( k_0 + \log x \big)}
\\
& = 2^{-3/2} \log \Big( 1 + k_0^{-1} \log \Big( 1 + \frac{s}{|a_n|} \Big) \Big)
\ge 2^{-3/2} \log \Big( 1 + k_0^{-1} \log \Big( 1 + \frac{s}{r} \Big) \Big).
\qedhere
\end{split}
\]
\end{proof}

\begin{proof}[Proof of Theorem~\ref{mainThm:qh}(2), for the quasihyperbolic metric]
We use the characterization of Bonk, Heinonen and Koskela \cite{BHK}.
Hence it suffices to show that the domain in not inner uniform.
So, suppose for a contradiction that the domain is $A$-inner uniform
for some fixed $A>0$.

We define $s_n:= \max_{1 \le m \le n} (b_m-a_m)$.
It is clear that $s_n$ is an increasing sequence
and $\lim_{n\to\infty}s_n/a_n=0$.
If we define $g_n:= \sqrt{s_n/a_n}\,$, then
$b_m-a_m \le a_ng_n^2$ for every $1 \le m \le n$
and $\lim_{n\to\infty} g_n=0$.

Since $g_n>0$, we can choose a subsequence $\{g_{n_k}\}$ with
$g_{n_k}\ge g_m$ for every $m \ge n_k$;
consider a fixed $n$ from the sequence $\{n_k\}$.
Set $c_n=\frac{b_n+a_n}2$, the mid-point of $(a_n,b_n)$.
We define $x_n=c_n + i c_n g_n$ and
$y_n=c_n - i c_n g_n$. Since $[x_n,y_n]\subset \Omega$,
we have $\ell_{\text{Eucl},\Omega}([x_n,y_n]) = 2 c_n g_n$.
Let $\gamma$ be an $A$-inner uniform curve joining
$x_n$ and $y_n$, and let $z\in \gamma\cap \R$.
Since $|x_n-z|, |y_n-z|\ge c_n g_n$, we conclude
by the uniformity of the curve that $\delta_\Omega(z) \ge \frac{c_n g_n}{A}$.
On the other hand, the uniformity of $\gamma$ also implies
that $|z-c_n| \le 2 A c_n g_n$.

We may assume that $n$ is so large that $c_n > 2 A c_n g_n$.
Then $z$ lies in the positive real axis, which means that
$z\in (a_m,b_m)$ for some $m\ge 1$. If $m\le n$, then we have
$b_m-a_m \le s_n = a_n g_n^2 < c_n g_n^2$. For $m>n$ we have
$b_m-a_m \le g_m^2 a_m \le g_n^2 a_m$. However,
since $a_m < z \le c_n + 2 A c_n g_n < 2 c_n$, so
for every $m$ we have $b_m-a_m < 2 c_n g_n^2$.

Since $\delta_\Omega(z) < \frac{b_m-a_m}2$, we conclude that
$\frac{c_n g_n}{A} < c_n g_n^2$. Since $g_n\to 0$ and
$A$ is a constant, this is a contradiction. Hence the
assumption that an $A$-inner uniform curve exists was false,
and we can conclude that the domain is not Gromov hyperbolic.
\end{proof}

For the proof in the hyperbolic case we need the following concepts.
A function between two metric spaces $f:X\longrightarrow Y$ is an
{\it $(a,b)$-quasi-isometry}, $a\ge 1$, $b\ge 0$, if
\[
\frac1a\,d_{X}(x_1,x_2)-b\le d_{Y}(f(x_1),f(x_2))
\le a d_{X}(x_1,x_2)+b\,,
\qquad \text{for every } x_1,x_2\in X .
\]
An $(a,b)$-\emph{quasigeodesic} in $X$ is an
$(a,b)$-quasi-isometry between an interval of $\R$ and $X$.

For future reference we record the following lemma:

\begin{lemma}
\label{l:quasigeodesic}
Let us consider a geodesic metric space $X$
and a geodesic $\g:I\longrightarrow X$, with $I$ any interval,
and $g:I\longrightarrow X$, with $d(g(t),\g(t))\le \e$ for every $t\in I$.
Then $g$ is a $(1,2\e)$-quasigeodesic.
\end{lemma}

\begin{proof}
We have for every $s,t\in I$
$$
d(g(s),g(t))
\ge d(\g(s),\g(t)) -d(\g(s),g(s)) -d(\g(t),g(t))
\ge  |t-s| -2\e .
$$
The upper bound is similar.
\end{proof}

\begin{proof}[Proof of Theorem~\ref{mainThm:h}(2), for the hyperbolic metric]
We consider two cases: either $\{b_m-a_m\}_m$ is bounded or unbounded.
We start with the latter case.

As in the previous proof, we define $s_n:= \max_{1 \le m \le n} (b_m-a_m)$
and $g_n:= \sqrt{s_n/a_n}\,$. Then $b_m-a_m \le a_ng_n^2$
for every $1 \le m \le n$ and $\lim_{n\to\infty} g_n=0$.
Since $g_n>0$, we can choose a subsequence $\{g_{n_k}\}$ with
$g_{n_k}\ge g_m$ for every $m \ge n_k$.
Since $\{b_m-a_m\}_m$ is not bounded we may, moreover,
choose the sequence so that $g_n^2 = (b_n-a_n)/a_n$
for every $n\in \{n_k\}$. Fix now $n$ from the sequence $\{n_k\}$.
As before, we conclude that
$b_m-a_m \le a_n g_n^2$ for $m\le n$ and
$b_m-a_m \le a_m g_m^2 \le a_m g_n^2$ for $m>n$.

\begin{figure}[ht]
\hspace{-0.8cm}
\includegraphics[scale=1.4]{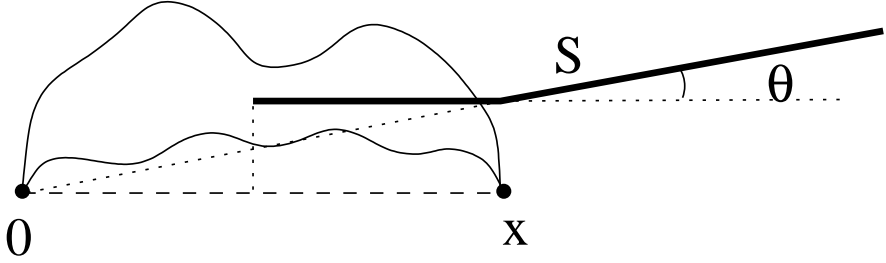}
\hspace{-1.5cm}
\caption{The set $S$}\label{figureS}
\hspace{-1.5cm}
\end{figure}

Consider $x\in (a_n,b_n)$ which lies on the shortest fundamental
geodesic $\gamma_n$ joining $(-\infty,0)$ with $(a_n,b_n)$.
Define an angle $\theta= \atan {g_n} \in (0,\pi/2)$ and a set
\[
S = [\tfrac12 x + i x{g_n}, x + i x{g_n}]
\cup \{x +  ix{g_n} + te^{\pi i\theta}\; |\; t \ge 0 \}.
\]
The set $S$ is shown in Figure~\ref{figureS}. Notice that any point
$\zeta\in S$ satisfies ${g_n} \Real \zeta \le \Imaginary \zeta \le 2{g_n} \Real
\zeta$. It is clear
that $\gamma_n$ hits the set $S\cup [\tfrac12 x + i x {g_n}, \tfrac12 x]$.
We claim that it in fact hits $S$. Assume to the contrary that this is not
the case. Then it hits $[\tfrac12 x + i x{g_n}, \tfrac12 x]$. Let $\gamma'$ denote a
part of $\gamma_n$ connecting $x$ and this segment which does not
intersect $S$. Since $\Omega$ is a Denjoy domain, we conclude that
$b\mapsto \lambda_\Omega(a + ib)$ is decreasing for $b>0$
(see \cite[Theorem 4.1(i)]{M}).
Hence $\ell_{h,\Omega}(\gamma') \ge \ell_{h,\Omega}([\tfrac12 x + i x{g_n}, x + i x
{g_n}])$. Since the gap size in $[\tfrac12 x, x]$ is at most
$a_n g_n^2$, we have $\delta_\Omega(w) \le \sqrt{x^2 {g_n}^2 + a_n^2 g_n^4}\le
\sqrt{2}\, x {g_n}$. Since the gap size is smaller than the distance
to the boundary, it follows from Theorem~\ref{BPthm} that
\[
\lambda_\Omega(w)
\ge \frac{C}{\delta_\Omega(w)}
\ge \frac{C}{x{g_n}}
\]
for $w\in [\tfrac12 x + i x{g_n}, x + i x {g_n}]$. Multiplying this
with the Euclidean length $\tfrac12 x $ of the segment gives
\[
\ell_{h,\Omega}(\gamma_n)
\ge \ell_{h,\Omega}([\tfrac12 x  + i x{g_n}, x + i x {g_n}])
\ge \frac{C}{ g_n}.
\]

We next construct another path $\sigma$ and show that it is in
the same homotopy class as the supposed geodesic, only shorter. Let
$z$ be the midpoint of gap $n$ and let $\sigma$ be the curve
$[z,z+iz]\cup[z+iz,-z+iz]\cup[-z+iz, -z]$.
Using $b_{n}-a_{n} = a_n g_n^2$ we easily calculate
\[
\ell_{h,\Omega}(\sigma)\le 2\ell_{k,\Omega}(\sigma)\le
2\log\Big(\frac{2z}{a_n g_n^2}\Big) + C \le 4\log\Big(\frac{1}{g_n}\Big) + C
\]
with an absolute constant $C$. The curve $\sigma$ joins $(-\infty,0)$ and
$(a_n,b_n)$; therefore $\ell_{h,\Omega}(\gamma_n)\le
\ell_{h,\Omega}(\sigma)$. But this contradicts the previously derived
bounds for the lengths as $g_n\to 0$.

Therefore the supposition that $\gamma_n$ does not
intersect $S$ was wrong, so we conclude that $\gamma_n \cap S\not=
\emptyset$. Let now $\zeta\in S\cap \gamma_n$. We claim that
$h_{\Omega}(\zeta,\R) \to \infty$,
which means the domain is not Gromov hyperbolic, by Theorem~\ref{th:2.1}.
Let $\xi\in \Omega\cap \R$; chose $m$ so that $\xi \in (a_m,b_m)$.
Let $\alpha$ be a curve joining $\xi$ and $\zeta$.

If $0<m\le n$, then the size of $(a_m,b_m)$ is at most $a_n g_n^2$, so $\delta_\Omega(\xi) \le
a_n g_n^2$. Then $\alpha$ has Euclidean length at least $\Imaginary \zeta \ge
x {g_n}$, so by Lemma~\ref{minLenLem},
$\ell_{h,\Omega}(\alpha)\ge c\log\log (C/g_n)$.
As $g_n\to 0$, this bound tends to $\infty$.
If, on the other hand, $m>n$,
then the Euclidean length of $\alpha$ is at least
\[
d(\xi,\zeta) \ge d(\xi,S) \ge \xi \sin \theta \ge \tfrac 12 \, \xi
\tan \theta= \tfrac12\, \xi g_n,
\]
and the size of the gap is at most $a_m g_n^2$.
By Lemma~\ref{minLenLem} this implies that
$\ell_{h,\Omega}(\alpha)\ge c\log\log (C/g_n)$.
As $g_n\to 0$, this bound again tends to $\infty$.

It remains to consider $m=0$, i.e., $\xi<0$.
We consider only the case $\zeta \in [\tfrac12 x + i x{g_n},
x + i x {g_n}]$, since the other case is similar.
Now the Euclidean length of $\alpha$ is at least $\tfrac12 x$.
Since the gap size in $[0, \tfrac12 x]$ is at most
$a_n g_n^2$, we see that the boundary satisfies the separation
condition when $|\Im z|\ge a_n g_n^2$
in which case also $\delta_\Omega(z)\ge |\Im z|\ge a_n g_n^2$.
Since $\lambda_\Omega(z)$ is decreasing in $|\Im z|$
(see \cite[Theorem~4.1(i)]{M}), we conclude that
\begin{equation}\label{negAxisEq}
\lambda_\Omega(z) \ge \frac{C}{\max\{|\Im z|, a_n g_n^2\}}
\ge \frac{C}{\max\{\delta_\Omega(z), a_n g_n^2\}}
\end{equation}
for the points on the curve with $\Re z \in (0,x/2)$.
Let $\alpha^-$ be the part of $\alpha$ on which
$\delta_\Omega(z)< a_n g_n^2$. If $\ell_{\text{Eucl}}(\alpha^-)> x g_n^{3/2}$, then
\[
\ell_{h,\Omega}(\alpha) \ge \ell_{h,\Omega}(\alpha^-)
\ge \frac{x g_n^{3/2}}{a_n g_n^2} > g_n^{-1/2}.
\]
If $\ell_{\text{Eucl}}(\alpha^-)\le x g_n^{3/2}$, then
$\ell_{\text{Eucl}}(\alpha\setminus \alpha^-) > \frac12x - x g_n^{3/2}$.
Hence we conclude (as in the proof of part $(1)$ in Lemma~\ref{minLenLem}) that
\[
\int_\alpha \lambda_\Omega(z)\, |dz| \ge C
\int_{\delta_\Omega(\zeta) + x g_n^{3/2}}^{x/2} \frac{dt}{t}
\ge C \log\Big(\frac{x/2}{\sqrt{2}\, a_n g_n + x g_n^{3/2}}\Big)
\ge C \log\Big(\frac{1}{g_n}\Big) - C.
\]
Hence in either case we get a lower bound which
tends to infinity as $g_n\to 0$.

\bigskip

This takes care of the case when $\{b_m-a_m\}_m$ is unbounded.
Assume next that $\sup_m (b_m-a_m) = M < \infty$. In this case it is difficult
to work with bigons, since we do not get a good control on what the
gedesics look like; the problem with the previous argument is
that we cannot choose $g_{n_k}^2= (b_{n_k}-a_{n_k})/a_{n_k}$ in our sequence,
and consequently we do not get a good bound on the length of the
curve $\sigma$, as defined above.

To get around this we consider a geodesic triangle.
Assume for a contradiction that $\hOmega$ is $\delta$-Gromov hyperbolic.
By geodesic stability \cite{BS}, there exists a number $\delta'$ so
that every $(\sqrt{2},0)$-quasigeodesic triangle is $\delta'$-thin.

Fix $R\gg M^2$ and set $w_\pm = \pm i R$. Let $\gamma_0$ be the geodesic
segment joining $w_+$ and $w_-$. Choose $t>0$ so large that
$h_{\Omega}(\gamma_0, H_t) > \delta'$, where $H_t=\{z\in \C\,|\, \Re z>t\}$.
Let $w\in \Omega\cap\R$ be a point in $H_{2\max\{t,R\}}$, and
let $\gamma_+\subset \overline{\H^2}$ be a geodesic joining
$w$ and $w_+$.

If $\gamma_+$ dips below the ray from $w$ through $w_+$,
then we replace the part below the ray by a part of the ray.
The resulting curve is denoted by $\tilde \gamma_+$.
Let us show that $\tilde \gamma_+$ is a quasigeodesic.
We define a mapping $f\colon \gamma_+ \to \tilde \gamma_+$
as follows. If $x\in \gamma_+ \cap \tilde \gamma_+$,
then $f(x)=x$. If $x\in \gamma_+ \setminus \tilde \gamma_+$
then we set $f(x)$ to equal the point on $\tilde \gamma_+$ with
real part equal to $\Re x$.

Since $\Omega$ is a Denjoy domain, the function
$b\mapsto \lambda_\Omega(a + ib)$ is decreasing for $b>0$
(see \cite[Theorem~4.1(i)]{M}). Hence $\lambda_\Omega(f(x)) \le
\lambda_\Omega(x)$. The arc-length distance element
is the vertical projection of the distance element at $x$ to
the line through $w$ and $w_+$: specifically, the
distance element $(dx,dy)$ becomes $(dx, \theta dx)$,
where $\theta$ is the slope of the line. Thus the maximal increase
in the distance element is $\sqrt{1+\theta^2}$. Since
the slope of the line lies in the range $[-1,0)$, we conclude from
these facts that $\tilde\gamma_+$ is a $(\sqrt{2},0)$-quasigeodesic.

Similarly, we construct $\tilde\gamma_-$ and conclude that
it is a $(\sqrt{2},0)$-quasigeodesic. Choose now
$\zeta\in \tilde\gamma_+ \cap H_{\max\{t,R\}}$ with $\Im \zeta = \sqrt{R}$.
Since $\gamma_0\cup \tilde\gamma_+ \cup \tilde\gamma_-$
is a $(\sqrt{2},0)$-quasigeodesic triangle,
it should be possible to to connect $\zeta$ with some
point in $\gamma_0\cup \tilde\gamma_-$ using a path of
length $\delta'$. By the definition of $t$, $h_\Omega(\zeta,\gamma_0)>\delta'$.
If $\alpha$ is a path connecting $\zeta$ and $\gamma_-$, then
it crosses the real axis at some point $\xi$.
If $\xi$ lies in $(a_m,b_m)$, $m>0$, then
$\ell_{h,\Omega}(\alpha)\ge C \log \log \frac{\sqrt{R}}{M}$, by Lemma~\ref{minLenLem}.
Otherwise, $\xi\in (-\infty,0)$. This case is handled as
in the first case of the proof, see the paragraph around
\eqref{negAxisEq}. In each case we see that
$h_\Omega(\zeta,\gamma_-)>\delta'$ provided $R$ is large enough.
But this means that $\Omega$ is not Gromow hyperbolic,
as was to be shown.
\end{proof}

In Theorem~\ref{mainThm:h}(2) the gaps $(a_n,b_n)$ and $(a_{n+1},b_{n+1})$
are separated by a boundary component $[b_n,a_{n+1}]$.
We easily see from the proofs that it would have made no
difference if this boundary component had some gaps, as long
as they at most comparable to the lengths of the adjecent gaps,
$(a_n,b_n)$ and $(a_{n+1},b_{n+1})$.
Thus we get the following stronger theorem by the same proofs.
(In the proofs we can assume that $(-\infty,0)\subset \Omega$,
by using Theorem~\ref{th:negativehalfaxis}).

\begin{theorem}
\label{HypQuasihyp}
Let $\Omega$ be a Denjoy domain with $\Omega\cap \R=\bigcup (a_n,b_n)$
and $\limsup_{n\to\infty}a_n=\infty$.
Suppose $G\colon \R^+ \to \R^+$ is a function with $\lim_{x\to\infty}G(x)=0$.
If $b_n-a_n \le a_n G(a_n)$ for every $a_n>0$, then $\hkOmega$, the hyperbolic or
quasihyperbolic metric, is not Gromov hyperbolic.
\end{theorem}

The function $G$ plays the role of $g_n^2$ in the proofs of
Theorem~\ref{mainThm:h}(2).

\begin{remark}
The condition $\Omega\cap\R = \bigcup (a_n,b_n)$
(without the hypothesis $b_n \le a_{n+1}$ for every $n$)
allows any topological behaviour; for instance,
$\partial\Omega$ can contain a countable
sequence of Cantor sets.
\end{remark}

Let $E_0 \subset [0,t)$ be closed, $t>0$, set $E_n:= E_0+tn$ for
$n\in \N$, and $\Omega:=\C \setminus \cup_{n} E_n$.
Then $\Omega$ satisfies the hypotheses of Theorem~\ref{HypQuasihyp}
for $G(x) = t/x$.
From this we deduce Corollary~\ref{periodicCor},
the non-hyperbolicity of periodic Denjoy domain,
in the case the index set is $\N$.
The case with index set $\Z$ follows from this
and Theorem~\ref{th:negativehalfaxis}.

%%%%%%%%%%%%%%%%%%%%%%%%%%%%%%%%%%%%%%%%%%%%%%%%%%%%%%%%%%%%%%%%%%%%%%%%%
%%%%%%%%%%%%%%%%%%%%%%%%%%%%%%%%%%%%%%%%%%%%%%%%%%%%%%%%%%%%%%%%%%%%%%%%%
%%%%%%%%%%%%%%%%%%%%%%%%%%%%%%%%%%%%%%%%%%%%%%%%%%%%%%%%%%%%%%%%%%%%%%%%%

\section{On the far side of the accumulation point}

\begin{lemma}
\label{l:endpoints}
Let $\Omega$ be a Denjoy domain with
$\Omega\cap\R = \cup_{n=0}^\infty (a_n,b_n)$ and $a_0=-\infty$.
If $\hOmega$ is not Gromov hyperbolic,
then for every $N>0$
there exist fundamental geodesics $\g_{n_k}$, $n_k > N$,
such that the hyperbolic distance of the endpoints
of $\g_{n_k}$ to $(-\infty,b_0)$ is greater than $N$,
and points $z_k\in \g_{n_k}$ with
$\lim_{k\to\infty} h_\Omega(z_k,\R)= \infty$.
\end{lemma}

\begin{proof}
Let us choose fundamental geodesics $\{\g_{n}^0\}$.
Since $\hOmega$ is not Gromov hyperbolic, by
Theorem~\ref{th:2.1}
there exists points  $w_k\in \g_{n_k}^0$ with $n_k > N$ and
$\lim_{k\to\infty} h_\Omega(w_k,\R)= \infty$. Since $\lim_{x\to b_n}
h_\Omega(x,(-\infty,b_0))= \infty$ for every $n$,
there exist $x_0 \in (a_{0},b_{0})$
and $x_{n_k} \in (a_{n_k},b_{n_k})$,
with $h_\Omega(x_0,(-\infty,b_0)),h_\Omega(x_{n_k},(-\infty,b_0))>N$.

Let us consider the fundamental geodesics $\g_{n_k}$ joining $x_0$
and $x_{n_k}$, as well as the bordered Riemann surface
$X:=\Omega \cap \overline{\H^2}$, which as in the proof of Theorem~\ref{l:halfplane}
can be shown to have $\log\big(1+\sqrt{2}\,\big)$-thin triangles.

Let $Q_k$ be the geodesic quadrilateral given by
$\g_{n_k}^0$, $\g_{n_k}$ and the two geodesics
(contained in $(a_{0},b_{0})$ and $(a_{n_k},b_{n_k})$)
joining their endpoints.
Since $Q_k \subset X$, it is $2\log\big(1+\sqrt{2}\,\big)$-thin,
and there exists $z_k \in \g_{n_k}\cup \R$
with $h_\Omega(z_k,w_k)\le 2\log\big(1+\sqrt{2}\,\big)$.

Since $\lim_{k\to\infty} h_\Omega(w_k,\R)= \infty$,
we deduce that $z_k \in \g_{n_k}$ for every $k\ge k_0$ and
$\lim_{k\to\infty} h_\Omega(z_k,\R)= \infty$.
\end{proof}

\begin{lemma}[Lemma 3.1, \cite{APR}]
\label{aprlemma}
Consider an open Riemann surface $S$ of hyperbolic type,
a closed non-empty subset $C$ of $S$, and set $S^*:=S\setminus C$.
For $\epsilon>0$ we have
$1 < \ell_{S^*}(\g)/\ell_S(\g) < \coth (\e/2)$,
for every curve $\g\subset S$ with finite length in $S$ such that
$h_S(\g,C)\ge\e$.
\end{lemma}

%\begin{definition}
Given a Riemann surface $S$, a geodesic $\g$ in $S$,
and a continuous unit vector field $\xi$ along $\g$ orthogonal to $\g$,
we define \emph{Fermi coordinates} based on $\g$ as the map
$Y(r,t):=\exp_{\g(r)} t \xi(r)$.
%\end{definition}

It is well known that if the curvature is $K\equiv -1$, then
the Riemannian metric can be expressed
in Fermi coordinates as $ds^2= dt^2 + \cosh\!^2 t \, dr^2$
(see e.g. \cite[p.~247--248]{C}).

\begin{corollary}
\label{c:geodesic}
Consider an open Riemann surface of hyperbolic type $S$,
a closed non-empty subset $C$ of $S$, and set $S^*:=S\setminus C$.
For $\epsilon>0$ and $C_\e:=\{z\in S: \, h_S(z,C) \ge \e \}$
we have
$$
\begin{aligned}
h_S(z,w) & \le h_{S^*}(z,w) , \qquad \text{ for every } \, z,w\in S^* ,
\\
h_{S^*}(z,w) & \le \coth (\e/2) \, h_{S|C_\e}(z,w) ,
\qquad \text{ for every } \, z,w\in C_\e \,.
\end{aligned}
$$
Furthermore, if $S$ is a Denjoy domain and
$C$ is a component of $S\cap \R$ then
$$
h_{S^*}(z,w) \le \cosh \e \, \coth (\e/2) \, h_{S}(z,w) ,
$$
for every $z,w$ in the same component of $C_\e$
with $\Imaginary z, \Imaginary w \ge 0$.
\end{corollary}

\begin{proof}
The first and second inequalities are direct consequences of Lemma
\ref{aprlemma}. In order to prove the third one, it is sufficient to prove that
\begin{equation}
h_{S|C_\e}(z,w) \le (\cosh \e) \, h_{S}(z,w),
\label{2.1}
\end{equation}
for every $z,w$ in the same component of $C_\e$
with $\Imaginary z, \Imaginary w \ge 0$.

Fix $z,w$ in the same component $\Gamma$ of $C_\e$.
Since $\Imaginary z, \Imaginary w \ge 0$ there exists
a unique geodesic $\g\subset S\cap \overline{\H^2}$
joining $z$ with $w$.

If $\g\subset \Gamma$, then $h_{S|C_\e}(z,w) = h_{S}(z,w)$. If $\g$ is not
contained in $\Gamma$, then it is sufficient to show that there exists a
curve $\eta$ joining $z$ and $w$ in $\Gamma$, with
$\ell_{h,S}(\eta) \le (\cosh \e) \, \ell_{h,S}(\g)$.
In order to prove this,
consider the geodesics $\g_z,\g_w\subset S\cap \overline{\H^2}$
joining $z$ and $w$ with $C$,
and the geodesic $\g_0\subset C$
joining the endpoints of $\g_z,\g_w$ (which are in $C$).

We denote by $P$ the simply connected closed region with boundary
$\g \cap \g_z \cap \g_w \cap \g_0$.
Since $P$ is simply connected, we can identify it with a domain
$P_0\subset \overline{\H^2}$ using Fermi coordinates based on $C$.

If $g$ is the lift of $\g$, then
$g_1:=g\cap \{(r,t):\,0\le t\le \e\}$ is the lift of $\g\setminus C_\e$.
If $g\cap \{(r,t):\, t= \e\}=\{(r_1,\e),(r_2,\e)\}$ (with $r_1<r_2$),
then we define $g_2:=\{(r,\e):\, r_1\le r \le r_2\}$ and
$g_0:=\{(r,0):\, r_1\le r \le r_2\}$.
Notice that in order to prove
(\ref{2.1})
it is sufficient to show that
$\ell(g_2) \le (\cosh \e) \, \ell(g_1)$.
But this is a direct consequence of the facts
$\ell(g_0) \le \ell(g_1)$ and $\ell(g_2) = (\cosh \e) \, \ell(g_0)$.
\end{proof}

\begin{proof}[Proof of Theorem~\ref{th:negativehalfaxis}]
Since $\hkOmega$ is not Gromov hyperbolic,
by Proposition~\ref{finitetype}, we conclude that
$\Omega$ has countably infinitely many boundary components:
$\Omega\cap\R = \cup_{n=0}^\infty (a_n,b_n)$. Without loss of generality we
can assume that $(-\infty,0)\subseteq (a_1,b_1)$.

We first prove that
$(\Omega\setminus F, k_{\Omega\setminus F})$ is not Gromov hyperbolic.
Let us consider fundamental geodesics $\g_n$ of
$k_{\Omega}$ joining the midpoint $c_0$ of $(a_0,b_0)$
with the midpoint $c_n$ of $(a_n,b_n)$ for $n\ge 2$
which are shortest possible.
Since $\g_n$ is contained in $\{z\in\C:\,c_0\le \Real z \le c_n\}$,
and $k_{\Omega\setminus F}=k_{\Omega}$ in
$\{z\in\C:\,\Real z \ge \inf_{n\ge 2} a_n \}$,
we deduce that $\g_n$ is also a fundamental geodesic with the metric
$k_{\Omega\setminus F}$.

Since $\kOmega$ is not Gromov hyperbolic,
there exist points $z_k\in \g_{n_k}$ with
$\lim_{k\to\infty} k_\Omega(z_k,\R)=\infty$ by
Theorem~\ref{th:2.1}. Since $\g_{n_k}$ are also fundamental geodesics
with the metric $k_{\Omega\setminus F}$, we deduce that
$\lim_{k\to\infty} k_{\Omega\setminus F}(z_k,\R)\ge \lim_{k\to\infty}
k_\Omega(z_k,\R)=\infty$. Consequently, $(\Omega\setminus F, k_{\Omega\setminus F})$ is not Gromov hyperbolic.
\smallskip

We now prove that
$(\Omega\setminus F, h_{\Omega\setminus F})$ is not Gromov hyperbolic.
Choose $\e_0>0$. Since $\hOmega$ is not Gromov hyperbolic,
by Lemma~\ref{l:endpoints} there exist fundamental geodesics $\g_{n_k}$
of $h_{\Omega}$,
such that the hyperbolic distance of the endpoints of $\g_{n_k}$ to
$(-\infty,b_1)$ is greater than $\e_0$, and points $z_k\in \g_{n_k}$ with
$\lim_{k\to\infty} h_\Omega(z_k,\R)= \infty$.

Fix $\e \in \big(0, \min\{\e_0, \min_k h_\Omega(z_k,\R)\}\big)$.
If we define
$$
U_\e:=\{z\in\Omega: \, h_\Omega(z,(-\infty,b_1))\ge \e \},
$$
we see that $z_k \in \g_{n_k}\cap U_\e$ for every $k$.
(Notice that $\g_{n_k}\cap \partial U_\e$ has at most two points.)
If $\g_{n_k}\cap \p U_\e$ is empty or a one-point set, we define $g_{n_k}:=\g_{n_k}$.
Since the endpoints of $\g_{n_k}$ are in $U_\e$,
we conclude that $g_{n_k}\subset U_\e$.

Then assume that $\g_{n_k}\cap \p U_\e=\{w^1,w^2\}$. If there is an arc
$\a$ in $\p U_\e$ joining $w^1$ and $w^2$, we define a curve $g_{n_k}$
joining $(a_0,b_0)$ with
$(a_{n_k},b_{n_k})$ in $U_\e$, by $g_{n_k}:=(\g_{n_k}\cap U_\e)\cup \a$.
Then $\g_{n_k}$ and $g_{n_k}$ have the same endpoints and are homotopic.
If there is not an arc in $\p U_\e$ joining $w^1$ and $w^2$,
there are still maximal arcs $\a,\b$ in $\p U_\e$ joining
$w^1$ and $\o^1 \in (a_{m^1},b_{m^1})$, and $w^2$ and $\o^2 \in
(a_{m^2},b_{m^2})$, respectively, and a geodesic $\eta$ (with respect to
$h_\Omega$) in $\Omega\setminus U_\e$ joining $\o^1$ and $\o^2$, such that
if $\g_{n_k} \cap U_\e=[z^1,w^1]\cup [z^2,w^2]$,
then $[z^1,w^1] \cup \a \cup \eta \cup \b \cup [z^2,w^2]$
has the same endpoints as $\g_{n_k}$, and they are homotopic.

Since $\e< h_\Omega(z_k,\R)$, we have either $z_k \in [z^1,w^1]$ or
$z_k\in [z^2,w^2]$. Without loss of generality we can assume that
$z_k\in [z^2,w^2]$. Then we define $g_{n_k}:=\b\cup [z^2,w^2]\subset U_\e$,
which is a curve joining $(a_{m^2},b_{m^2})$ with $(a_{n_k},b_{n_k})$.

In any case, Lemma~\ref{l:quasigeodesic} gives that $g_{n_k}$ is a
$(1,2\e)$-quasigeodesic with respect to $h_\Omega$. Hence, for every
$t,s,$ we have
$$|t-s| -2\e \le h_{\Omega}\big( g_{n_k}(t), g_{n_k}(s) \big)
\le |t-s| +2\e .$$
Since $g_{n_k}$ is contained in $U_\e$,
Corollary~\ref{c:geodesic} implies that
$$
\begin{aligned}
|t-s| -2\e
& \le h_{\Omega}\big( g_{n_k}(t), g_{n_k}(s) \big)
< h_{\Omega\setminus F}\big( g_{n_k}(t), g_{n_k}(s) \big)
\\
& \le h_{\Omega\setminus (-\infty,0]}\big( g_{n_k}(t), g_{n_k}(s) \big)
\\
& \le \cosh \e \, \coth (\e/2) \,h_{\Omega}\big( g_{n_k}(t), g_{n_k}(s) \big)
\\
& \le \cosh \e \, \coth (\e/2) \,\big(|t-s| +2\e \big) ,
\end{aligned}
$$
and hence $g_{n_k}$ is a $\big( \cosh \e \coth(\e/2), 2 \e
\cosh \e \coth(\e/2) \big)$-quasigeodesic with respect to $h_{\Omega\setminus F}$.
% Since $\e< h_\Omega(z_k,\R)$ and $g_{n_k}$ is contained in $U_\e$,
% we have $z_k \in g_{n_k}$ for every $k$.

To get  a contradiction, assume that $(\Omega\setminus F,
h_{\Omega\setminus F})$ is Gromov hyperbolic. Consider the fundamental
geodesic $\eta_{n_k}$ of $h_{\Omega\setminus F}$ with the same
endpoints as $g_{n_k}$. Then there is a constant $C$ such that
the Hausdorff distance of $g_{n_k}$ and $\eta_{n_k}$ is less than $C$.
Hence, there exist points $w_k \in \eta_{n_k}$
with $h_{\Omega\setminus F}(z_k,w_k)\le C$, and thus
$$
\lim_{k\to\infty} h_{\Omega\setminus F}(w_k,\R)
\ge \lim_{k\to\infty} h_{\Omega\setminus F}(z_k,\R)- C
\ge \lim_{k\to\infty} h_{\Omega}(z_k,\R)-C= \infty ,
$$
which contradicts
$h_{\Omega\setminus F}$ being Gromov hyperbolic.
\end{proof}

%%%%%%%%%%%%%%%%%%%%%%%%%%%%%%%%%%%%%%%%%%%%%%%%%%%%%%%%%%%%%%%%%%%%%%%%%
%%%%%%%%%%%%%%%%%%%%%%%%%%%%%%%%%%%%%%%%%%%%%%%%%%%%%%%%%%%%%%%%%%%%%%%%%
%%%%%%%%%%%%%%%%%%%%%%%%%%%%%%%%%%%%%%%%%%%%%%%%%%%%%%%%%%%%%%%%%%%%%%%%%

%\bibliographystyle{amsplain}

\end{document}